\documentstyle{article}

\begin{document}

\centerline {\bf On (3,3)-homogeneous Greechie Orthomodular
Posets}
\vskip 0.3 cm
\centerline { Foat ~F.~ Sultanbekov
\footnote{ Department of Mathematics and Mechanics, Kazan State University,
Kremlevskaya Street, 18. Kazan, Tatarstan, Russia, 420008.}}

\vskip 0.5 cm
\noindent
{\bf ABSTRACT.} We describe (3,3)-homogeneous orthomodular posets for
some cardinality of their sets of atoms. We examine a state
space and a set of two-valued states of such logics. Particular
homogeneous OMPs with exactly $k$ pure states ($k=1,...,7, 10,11$)
have been constructed.
\vskip 0.5 cm

\noindent
{\bf 1. INTRODUCTION}
\vskip 0.2 cm

Homogeneous orthomodular posets (OMPs) are important [10], [13].
They can be used in constructing counterexamples or OMPs with
certain properties of the state space or the automorphisms group
[5], [6], [7].

Let $n, m$ be natural numbers. An OMP, $L$, is called {\it $(n,
m)$-homogeneous} ($(n, m)$-hom.), if its every atom is
contained in $n$ maximal, with respect to inclusion, orhogonal sets
of atoms (called {\it blocks}), and every such set of atoms
of $L$ is $m$-element. The well known concrete logics of the form
${\cal L}^{p}_q = \{ X \subset \{ 1, \ldots, pq\} \mid \hbox
{card} X \equiv 0 \, (\hbox{ mod} \, q)\}$ [10] are nice examples of
homogeneous OMPs. (3,3)-homogeneous logics arise when we
consider relational OMPs [4] on a finite set. Orthomodular
lattices of the kind were examined in [3],[12].

Let $L$ be a finite $(n, m)$-hom. OMP, $A$ the set of all atoms in $L$,
$B$ the set of all blocks in $L$, $S=S(L)$ the
set of all states on $L$, and $S_2$ the set of all two-valued states on
$L$. A state $s$ on $L$ is called {\it pure} if $s$ is an extreme
point of the convex set $S(L)$. It is easy to see that
$n \cdot \mbox{card}A = m \cdot \mbox{card}B $.

\vskip 0.2 cm
{\bf Theorem 1.1.}([10]){\it Suppose that $S_2 \not= \emptyset $ and
$f \in S_2$. Then} \mbox{card}$A =mk$, \mbox{card}$B=nk$ {\it where}
$k=$\mbox{card}$ (f^{-} (1) \cap A)$.
\vskip 0.2 cm

Let us recall some definitions of the theory of concrete logics
[8], [9]. Let $\Omega$ be a set and ${\cal P}(\Omega)$ the Boolean
algebra of all subsets of $\Omega$. A {\it concrete logic} (c.l.)
on $\Omega$ is a subset ${\cal E}$ of ${\cal P}(\Omega)$
satisfying $ (1)\,\, \Omega \in {\cal E}; \quad (2)\,\, x \in
{\cal E} \Rightarrow \Omega \backslash x \in {\cal E}; \quad
(3)\,\, x,y \in {\cal E}, x \cap y = \emptyset \Rightarrow x \cup
y \in {\cal E}.$

Denote by $V({\cal E})$ the real vector space of all signed
measures on ${\cal E}$ and put $ {\cal E}^{\circ } = \{ \mu \in
V({\cal P}(\Omega))| \forall x\in {\cal E} \,(\mu(x)=0)\}$ .
${\cal E}$ is called {\it regular} if every signed measure on
${\cal E}$ extends to a signed measure on ${\cal P}(\Omega)$.

\vskip 0.2 cm
{\bf Theorem 1.2.}([8], [9]){\it A concrete logic, ${\cal E}$, is
regular iff}  $\hbox{dim} {\cal E}^{\circ} + \hbox{dim} V({\cal E})
= \hbox{card}\,\Omega$.
\vskip 0.2 cm

A set $T \subset S_2(L)$ is called {\it full} if for all $ x,y
\in L $ it holds $ x \le y \Leftrightarrow \forall f \in T (f(x) \le f(y)$.
OMP $L$ is isomorphic to a c.l. iff $S_2(L)$ is full. In this
case, the sets $x'= \{ \mu \in S_2(L)| \mu (x) = 1 \}$ form a c.l.
${\cal E}={\cal E}(L)$ on $\Omega=S_2(L)$ wich is called a
{\it total representation} of $L$. We also mention a simple
criterion of the fullness of a $S_{2}$. The condition equivalent to
the fullness of $S_{2}$ is as follows: if $ x, y \in A$ are not
orthogonal, then there exists $s \in S_2(L)$ such that
$s(x)=s(y)=1$.

Let us dwell (3,3)-hom. finite OMPs, $L$. Let
$l_n =l(P_0,...,P_{n-1}, Q_0,...,Q_{n-1})$ be a loop [1] of order
$n$, where $P_i$ denote the atoms of $l_n$ lying in the
vertices of a $n$-polygon and $Q_i$ are the atoms lying in the
middles of sides of the $n$-polygon. So $\{ P_i, Q_i, P_{i+1} \}, i =
0,...,n-1$ (indices modulo $n$) are
all blocks of $l_n$. It is easy to see that a two-valued state,
$s$, on $l_n$ is well determined by $s(P_0),...,
s(P_{n-1})$. We use the following abbreviation: $P_{01}Q_0 =
\{P_0, P_1, Q_0 \}$.

\vskip 0.6 cm
\noindent
{\bf 2. (3,3)-HOMOGENEOUS FINITE OMPs}
\vskip 0.2 cm

Let $L$ be a (3,3)-hom. OMP. Obviously $\hbox{card}\,A =
\hbox{card}\,B \ge 15$ and $\hbox{card}\,A \not= 16$.
\vskip 0.1 cm

{\bf Theorem 2.1.} 1.{\it If} $\hbox{card}\,A = 15$ {\it then $L$
is isomorphic to ${\cal L}^{3}_2$.}

\noindent
2. {\it There exist} (3,3)-{\it hom. OMPs $L_{17}, L_{27}$ with}
\mbox {card}$A \in \{ 17, 27 \}$. {\it For $L_{17}$ the set $S_2 =
\emptyset$ and $S(L_{17})$ is isomorphic to a segment $[0;
\frac{2}{3}]$. For OMP $L_{27}$ the set $S_2$ is full and total
representation of $L_{27}$ is regular. }
\vskip 0.2 cm

\noindent
{\it Proof.} 1. Let \mbox {card}$A$ = 15. Using Greechie diagram
for $L$ it is easy write out all two-valued states of $L$. So,
\mbox {card}$S_2=6$ and $S_2$ is full. The total representation
of $L$ is minimal and isomorphic to ${\cal L}^{3}_2$.

\noindent
2. First we construct $L_{17}$. Let us consider loop $l_7=
l(P_0,...,P_6, Q_0,...,Q_6)$ and add the atoms $R_0, R_1, R_2;
R_i\not\in l_7 (i=0, 1, 2)$. For seven blocks of $l_7$ we add
following ten blocks:

\centerline { $P_{04}R_0,\, Q_{025},\, Q_{03}R_2,\, P_1Q_4R_1,\,
Q_{136}$,}
\centerline {$ Q_{15}R_0,\, P_{25}R_2, \, Q_{246}, \, P_{36}R_1, \,
R_{012}$.}
\vskip 0.1 cm

Next we prove that $S(L_{17})$ is isomorphic to a segment $[0;
\frac {2}{3}]$. Let $s \in S(L_{17})$. Put $s(P_0)=x,
s(Q_0)=y, s(Q_6)=z, s(R_0)=t$ and $s(Q_2)=u $. We show that
all values of $s$ are described by $x$.

\noindent
1) We have $s(P_1) = 1-x-y, s(P_6) = 1-x-z, s(P_4) = 1-x-t$ and
$s(Q_5) = 1-u-y, s(Q_4) = 1-u-z$. From the blocks
$ P_{45}Q_4, P_{56}Q_5 $ it follows that $s(P_5) =x+t+u+z-1$ and
$t=y$.

\noindent
2) Now $s(Q_1) = 1-s(Q_5)-s(R_0) = u, s(P_2) = x+y-u $ and $
s(R_1) =1-s(P_1)-s(Q_4) = x+y+z+u-1$. So $s(R_2)= 1-s(P_1)-s(Q_4)
=2-x-2y-u-z $. From the block $P_{25}R_2$ we get $u=x$.

\noindent
3) Next $s(R_1) = 1-s(P_1)-s(Q_4) = 2x+y+z-1$. From the blocks
$P_{36}R_1, Q_{136}$ we calculate $s(P_3)=1-x-y,s(Q_3)=1-x-z $.
Then from the block $P_{34}Q_3$ we have $1 =3-3x-2y-z$, or $ z=
2-2y-3x$.

\noindent
4) If $z$ from 3) is placed to  $s(R_2)$ and $s(Q_3)$ then we get
$s(R_2)=x, s(Q_3)=2x+2y-1$. So from the blocks
$Q_{03}R_2$ we have $x+y=\frac{2}{3}$.

So the state $s$ has only three values -- $x,\,\,
\frac{2}{3}-x$,\,\, and \, \, $\frac {1}{3}$, namely:

\vskip 0.1 cm
\centerline {$x$ -- on the atoms $P_0, Q_1, Q_2, R_2$; \,\,\,
$(\frac{2}{3}-x)$ -- on the atoms $P_2, Q_0, Q_6, R_0$}
\centerline {$\frac{1}{3}$ -- on all other remaining atoms.}
\vskip 0.2 cm

So $S(L_{17})$ is isomorphic to the segment  $[0; \frac{2}{3}]$
\vskip 0.1 cm

The Greechie diagram of the OMP $L_{27}$ is a cube in
three-dimensional space with 3 atoms on each edge, with one atom
in center of each side of the cube and with the last atom in the
center of cube. The blocks of $L_{27}$ are the lines drawing
parallel all axes throw the atoms. Thus, Greechie diagram of
$L_{27}$ is divided to three layers. Then state $s \in S_2 $ is
called {\it type 1 (type 2)} if $s$ equals 1 on main (secondary)
diagonal in one of the layers. Then $S_2$ has 6 states type 1 and
6 states type 2. So, \mbox {card}$S_2$=12. Next $\hbox{dim} {\cal
E}^{\circ} =3 $ and $ \hbox{dim} V({\cal E})=9$, where $ {\cal
E}$ is the total representation of $L_{27}$. By theorem 1.2. ${\cal
E}$ is regular.

\vskip 0.3 cm
\noindent
{\bf Theorem 2.2.} {\it There exist (3,3)-hom. OMPs with}
$\hbox{card}\,A \le 19$ {\it and with exactly $k$ pure states}
$(k=1,2,...,7,10,11)$.
\vskip 0.2 cm

\noindent
{\it Proof.} Let us denote by $H_{k}(m)$ a (3,3)-hom. logic with
$\hbox{card} A = m$ and $k$ pure states of $S$. Next, we
construct a nine OMPs: $H_{1}(19)$, $H_{2}(17)$, $H_{3}(18)$
$H_{4}(19)$, $H_{5}(19)$, $H_{6}(19)$, $H_{7}(18)$,
$H_{10}(18)$, $H_{11}(19)$.

We enumerate atoms of $H_{k}(m)$ by a natural numbers $1,
2,..., m$ and for a block $\{ i, j, n \}$ use abbreviation $i-j-n$.
Obviously every such OMP has the following 7 blocks: $B_{1}, ...,
B_{7}$: 1--2--3,\, 1--4--5,\, 1--6--7,\, 2--8--9,\, 2--10--11,\,
3--12--13,\, 3--14--15.
\vskip 0.1 cm

\noindent
{\bf 1) $H_{1}(19)$.} To $B_{1}, ..., B_{7}$ we add the following
12 blocks:

\vskip 0.1 cm
\centerline {4--8--12,\,\, 4--10--14,\,\, 5--9--16,\,\,
5--11--17,\,\, 6--8--15,\,\, 6--13--16,}
\centerline {7--9--18,\,\, 7--14--17,\,\, 10--16--19,\,\,
11--13--18,\,\, 12--17--19,\,\, 15--18--19.}
\vskip 0.1 cm

From the systen of linear equations $s(i)+s(j)+s(n)=1, \{ i, j, n
\} \in B$ we found the unique solution $s(i)=\frac {1}{3}
(i=1,..., 19)$.
\vskip 0.1 cm

\noindent
{\bf 2) $H_{2}(17)$.} Consider $L_{17}$ from Theorem 2.1 as
$H_{2}(17)$. Then $S$ has two pure states ($x=0, x= \frac
{2}{3}$).
\vskip 0.2 cm

\noindent
{\bf 3) $H_{3}(18)$.} To $B_{1}, ..., B_{7}$ we add the following
11 blocks:

\vskip 0.1 cm
\centerline {4--8--12, \,\, 4--14--16,\,\, 5--10--13,\,\,
5--17--18,\,\, 6--11--12,\,\, 6--16--18,}
\centerline {7--9--17,\,\, 7--10--15,\,\, 8--15--18,\,\,
9--13--16,\,\, 11--14--17.}
\vskip 0.1 cm

\noindent
Let $s \in S$. Put $s(17)=x, s(18)=y$. Then $s$ has values:
$x$ -- on the atoms 2, 4, 6, 13, 15; $y$ -- on the atoms
1, 9, 10, 12, 14; $1-x-y$ -- on the atoms 3, 5, 7, 8, 11, 16.
The state space $S$ is isomorphic to a triangle: $0 \le x \le 1,
0 \le y \le 1, x+y \le 1$. So, $S$ has 3 pure states:
(0,0),(1,0),(0,1).
\vskip 0.2 cm

\noindent
{\bf 4) $H_{4}(19)$.} To $B_{1}, ..., B_{7}$ we add the following
12 blocks:

\vskip 0.1 cm
\centerline {4--8--12,\,\, 4--10--14,\,\, 5--9--13,\,\,
5--11--16,\,\, 6--8--15,\,\, 6--11--17,}
\centerline {7--9--18,\,\, 7--10--19,\,\, 12--16--18,\,\,
13--17--19,\,\, 14--17--18,\,\, 15--16--19.}
\vskip 0.1 cm

\noindent
Let $s \in S$. Put $s(6)=x, s(8)=y$. Then $s$ has values:
$s(1)=2y-\frac {1}{3}, s(2)=1-2y, s(3)=s(12)=s(13)=\frac {1}{3}$,
$s(4)=s(5)= \frac {2}{3}-y, s(7)=\frac {4}{3}-x-2y, s(9)= y$,
$s(10)=\frac {2}{3}-x, s(11)=x+2y-\frac {2}{3},
s(14)=s(18)=x+y-\frac {1}{3}, s(15)=s(16)=1-x-y, s(17)=\frac
{5}{3}-2x-2y, s(19)=2x+2y-1$.

The state space $S$ is isomorphic to a parallelogram: $0 \le x
\le \frac {2}{3}, \frac {1}{6} \le y \le \frac {1}{2}, \frac {1}{2}
\le x+y \le \frac {5}{6}$. So, $S$ has 4 pure states:

\centerline {$(0,\frac {1}{2}),\,(\frac {1}{3},\frac {1}{6}),\, (\frac
{2}{3}, \frac {1}{6}),\, (\frac {1}{3}, \frac {1}{2})$.}
\vskip 0.2 cm

\noindent
{\bf 5) $H_{5}(19)$.} To $B_{1}, ..., B_{7}$ we add the following
12 blocks:

\centerline {4--8--12,\,\, 4--10--14,\,\, 5--9--15,\,\,
5--11--13,\,\, 6--13--16,\,\, 6--15--17,}
\centerline {7--11--18,\,\, 7--12--19,\,\, 8--17--18,\,\,
9--16--19,\,\, 10--17--19,\,\, 14--16--18.}
\vskip 0.1 cm

\noindent
Let $s \in S$. Put $s(15)=x, s(19)=y$. Then $s(18)=y, s(6)=\frac
{1}{3}-x+y $ and $s$ has values:

\centerline {$x$ -- on the atoms 1, 13;\,\, $(\frac {2}{3}-x)$ -- on the
atoms 3, 5;}
\centerline {$(\frac{2}{3}-y)$ -- on the atoms 7, 16, 17;\,\,
$\frac{1}{3}$ -- on all other remaining atoms.}

The state space $S$ is isomorphic to a pentagon: $0 \le x
\le \frac {2}{3}, 0 \le y \le \frac {2}{3}, x-y \le \frac {1}{3}$.
So, $S$ has 5 pure states: $(0, 0),\, (0, \frac {2}{3}),\,
(\frac {1}{3}, 0), \,(\frac {2}{3},\, \frac
{1}{3}),\,(\frac {2}{3}, \frac {2}{3})$.
\vskip 0.2 cm

\noindent
{\bf 6) $H_{6}(19)$.} To $B_{1}, ..., B_{7}$ we add the following
12 blocks:

\vskip 0.1 cm
\centerline {4--8--12,\,\, 4--10--14,\,\, 5--9--13,\,\,
5--11--16,\,\, 6--9--15,\,\, 6--11--17,}
\centerline {7--10--18,\,\, 7--13--19,\,\, 8--16--19,\,\,
12--17--18,\,\,
14--17--19,\,\, 15--16--18.}
\vskip 0.1 cm

\noindent
Let $s \in S$. Put $s(2)=x, s(10)=y$. Then $s(5)=x, s(6)=x-y+\frac
{1}{3}, s(11)=1-x-y, s(17)=2y-\frac {1}{3}, s(18)=1-2y$ and $s$
has also values:

\vskip 0.1 cm
\centerline {$y$ -- on the atoms 7, 15, 16;\,\, $(\frac {2}{3}-x)$ --
on the atoms 1, 9;}
\centerline {$(\frac{2}{3}-y)$ -- on the atoms 14, 19;\,\,
$\frac{1}{3}$ -- on the atoms 3, 4, 8, 12, 13.}
\vskip 0.1 cm

The state space $S$ is isomorphic to a hexagon: $0 \le x
\le \frac {2}{3}, \frac {1}{6} \le y \le \frac {1}{2}, y+x \le 1,
y-x \le \frac {1}{3} $. So, $S$ has 6 pure states:
\vskip 0.1 cm
\centerline {
$(0, \frac {1}{6}),\,\, (0, \frac {1}{3}),\,\,
(\frac {1}{6}, \frac {1}{2}),\,\, (\frac {1}{2}, \frac {1}{2}),\,\,
(\frac {2}{3}, \frac {1}{3}),\,\, (\frac {2}{3}, \frac {1}{6})$.}
\vskip 0.2 cm

\noindent
{\bf 7) $H_{7}(18)$.} To $B_{1}, ..., B_{7}$ we add the following
11 blocks:

\vskip 0.1 cm
\centerline {4--8--12,\,\, 4--10--14\,\, 5--11--15,\,\,
5--16--17,\,\, 6--8-18,\,\, 6--10--16,}
\centerline {7--9--15,\,\, 7--12--17,\,\, 9--13--16,\,\,
11--13--18,\,\, 14--17--18.}
\vskip 0.1 cm

\noindent
Let $s \in S$. Put $s(1)=x, s(3)=y, s(18)=z$. Then $s(2)=1-x-y,
s(17)=1-2z $ and $s$ has also values:
$x$ -- on the atoms 8, 10;\, $y$ -- on the atoms 9, 11;\,$(1-x-z)$ --
on the atoms 4, 6; \,$z$ -- on the atoms 5, 7, 12, 14, 16;\,
$(1-y-z)$ -- on the atoms 13, 15.
\vskip 0.1 cm

The state space $S$ is isomorphic to a polytope in
three-dimensional space: $0 \le x \le 1 , 0 \le y \le 1, 0 \le z
\le \frac {1}{2}, x+y \le 1, x+z \le 1, y+z \le 1 $. So, $S$ has
7 pure states:

\vskip 0.1 cm
\centerline {(0,0,0),\, (0, 1, 0),\,(1, 0, 0),\,
$(0, 0, \frac {1}{2}),\, (0, \frac {1}{2}, \frac {1}{2}),\,
(\frac {1}{2}, 0, \frac {1}{2}),\,
(\frac {1}{2}, \frac {1}{2}, \frac {1}{2})$.}
\vskip 0.1 cm

{\bf 8) $H_{10}(18)$.} To $B_{1}, ..., B_{7}$ we add the following
11 blocks:

\vskip 0.1 cm
\centerline {4--8--12,\,\, 4--10--14\,\, 5--8--16,\,\,
5--11--17,\,\, 6--12--16,\,\, 6--14--17,}
\centerline {7--9--13,\,\, 7--11--15,\,\, 9--17--18,\,\,
10--13--18,\,\, 15--16--18.}
\vskip 0.1 cm

\noindent
Let $s \in S$. Put $s(6)=x, s(16)=y, s(17)=z$. Then $s(1)= 2y-x,
s(11)=2y-z, s(12)=1-x-y, s(14)=1-x-z, s(11)=1-y-z $ and $s$
has also values:
\centerline {$x$ -- on the atoms 3, 4; $y$ --
on the atoms 8, 9, 13;$(1-2y)$ -- on the atoms 2, 5, 7;}
\centerline {$z$ -- on the atoms 5, 10;
$(1-y-z)$ -- on the atoms 13, 15.}
\vskip 0.1 cm

The state space $S$ is isomorphic to a polytope in
three-dimensional space: $0 \le x \le \frac {2}{3}, 0 \le y
\le \frac {1}{2}, 0 \le z \le \frac {2}{3}, x+y \le 1, x+z \le 1,
y+z \le 1, \frac {x}{2} \le y,  \frac {z}{2} \le y $.

So, $S$ has 10 pure states:

\vskip 0.1 cm
\centerline{
(0,0,0),\,
$(0, \frac {1}{2}, 0),\,
(\frac {1}{2}, \frac {1}{2}, 0),\,
(\frac {1}{2}, \frac {1}{2}, \frac {1}{2}),\,
(0, \frac {1}{2}, \frac {1}{2}),$ }
\vskip 0.2 cm

\centerline{
$(0, \frac {1}{3}, \frac {2}{3}),\,
(\frac {1}{3}, \frac {1}{3}, \frac {2}{3}),\,
(\frac {1}{2}, \frac {1}{4}, \frac {1}{2}),\,
(\frac {2}{3}, \frac {1}{3}, \frac {1}{3}),\,
(\frac {2}{3}, \frac {1}{3}, 0) $.}
\vskip 0.2 cm

{\bf 9) $H_{11}(19)$.} To $B_{1}, ..., B_{7}$ we add the following
12 blocks:

\vskip 0.1 cm
\centerline {4--8--12,\,\, 4--10--14\,\, 5--9--15,\,\,
5--11--13,\,\, 6--8--16,\,\, 6--10--17,}
\centerline {7--13--19,\,\, 7--15--18,\,\, 9--17--19,\,\,
11--16--18,\,\, 12--17--18,\,\, 14--16--19.}
\vskip 0.1 cm

\noindent
Let $s \in S$. Put $s(10)=x, s(15)=y, s(19)=z$. Then $ s(1)= x+y
-\frac {1}{3}, s(6)= \frac {1}{3}-x+z, s(7)=1-y-z $ and $s$
has also values:

\vskip 0.1 cm
\centerline {$x$ -- on the atom 8; $y$ --
on the atom 13;$z$ -- on the atom 18;}
\centerline {$\frac {2}{3}-x$ -- on
the atoms 2, 4; $\frac {2}{3}-y$ -- on the atoms 3, 5;
$\frac {2}{3}-z$ -- on the atoms 16, 17;}
\centerline {
$\frac {2}{3}-z$ -- on the atoms 16, 17;
$\frac {1}{3}$--on the atoms 9, 11, 12, 14. }
\vskip 0.1 cm

The state space $S$ is isomorphic to a polytope in
three-dimensional space: $0 \le x \le \frac {2}{3}, 0 \le y
\le \frac {2}{3}, 0 \le z \le \frac {2}{3}, x+y \ge \frac {1}{3},
z \ge x- \frac {1}{3}1, y+z \le 1 $.

So, $S$ has 11 pure states:

\vskip 0.1 cm
\centerline{
$(\frac {1}{3}, 0, 0),\,
(0, \frac {1}{3}, 0),\,
(0, \frac {2}{3}, 0),\,
(\frac {1}{3}, \frac {2}{3}, 0),\,
(\frac {2}{3}, \frac {2}{3}, \frac {1}{3}),\,
(\frac {2}{3}, 0, \frac {1}{3}),$ }
\vskip 0.2 cm

\centerline{
$(0, \frac {2}{3}, \frac {1}{3}),\,
(\frac {1}{3}, 0, \frac {2}{3}),\,
(\frac {2}{3}, 0, \frac {2}{3}),\,
(\frac {2}{3}, \frac {1}{3}, \frac {2}{3}),\,
(0, \frac {1}{3}, \frac {2}{3}) $.}

\vskip 0.3 cm
\noindent
{\bf Remark 2.3.} The special interest have (3,3)-hom. OMPs with a
unique state. What least number of atoms of such logic? An example with
22 atoms till now was known [2]. The example, constructed by us, has 19 atoms,
and we did not manage to construct OMP with smaller number of atoms.
Probably number of atoms 19 cannot be reduced; we yet have no the proof it.
In [11] was is developed a method of construction OMPs with a unique state,
however these of logics are not homogeneous.

Using Greechie loop lemma [1] it is not difficult to show, that all
listed above (3,3)-hom. logic $H_{k}(m)$ are not
orthomodular lattices. Certainly the examples of orthomodular lattices with
such property will have the much greater number of atoms.

\vskip 0.1 cm
\noindent
{\bf Remark 2.4.} There is a well-known method of constructing
the finite (3,3)-hom. OMPs with even card$A$. Let $n
\ge 9$ and $A = \{ a_{i}| i= 0,1, ..., 2n-1 \}$ be a set of
atoms. Then sets $\{ a_{2i}, a_{2i+1}, a_{2i+2} \};
\{ a_{2i-5}, a_{2i}, a_{2i+5} \} $ (indices modulo 2n) as the
blocks generate some (3,3)-hom. logic $L(2n)$. For example,
$L(22)$ has one, $L(20)$ has two, and $L(18)$ has three pure states
for the corresponding state spaces.
But the author is not familiar with any convenient method of constructing
such OMPs with odd card$A$.

\vskip 0.1 cm
\noindent
{\bf Remark 2.5.} (3,3)-hom. logics arise when we consider a
relational OMPs [4] on a finite set. For example, the relational
OMP on 8-element set is (3,3)-homogeneous. Every horizontal
summand of this OMP has 28 atoms and exactly one pure state.
We present all blocks of this summand:

\centerline {1--13--15,\, 1--6--21,\, 1--14--22,\,
2--6--25,\, 2--5--17,\, 2--18--26,\, 3--9--13,}
\centerline {3--10--25,\, 3--14--26,\, 4--9--17,\,
4--10--21,\, 4--18--22,\, 5--15--19,\, 6--23--27,}
\centerline {7--11--13,\, 7--12--25,\, 7--15--27,\,
8--11--21,\, 8--12--17,\, 8--19--23,\, 9--16--20,}
\centerline {10--24--28,\, 11--16--24,\, 12--20--28,\,
14--24--27,\, 15--20--26,\, 16--19--22,\, 18--23--28.}

\vskip 0.4 cm
\noindent
{\bf ACKNOWLEDGMENTS.}

This reseach was supported by RFBR (grant no. 01-01-00129) and
URFR (grant no. UR.04.01.061).

\vskip 0.6 cm
\noindent
{\bf REFERENCES}
\vskip 0.1 cm

1. Kalmbach, G., {\it Orthomodular Lattices}, Academic Press,
London, (1983).

2. Greechie, R.J, Miller, F.R., {\it On structures related to
states on an empirical logic I. Weights on finite spaces}, Technical
Report 16, Dept. of Math., Kansas State Univ., Manhattan, Kansas,
(1970), 1--25.

3. Kohler, E., {\it Orthomodulare Verbande mit
Regularitatsbedingungen (Orthomodular lattices with regularity
conditions)}, Journal of Geometry, Vol. 119 (1982), 129--145 (in
German).

4. Harding, J. {\it Decompositions in quantum logic},
Trans. Amer.Math. Soc., 348 (1996), 1839--1862.

5. Navara, M., and Rogalewicz, V., {\it Constructions of
orthomodular lattices with given state spaces}, Demostratio
Math., 21 (1988), 481--493.

6. Navara, M., and Tkadlec, J., {\it Automorphisms of
concrete logics}, Comment. Math. Univ. Carolinae
Math., 32 (1991), 15--25.

7. Navara, M., {\it An orthomodular lattice admiting no
group-valued measure}, Proc. Amer.Math. Soc.,122(1994), 7--12.

8. Ovchinnikov, P.G. and Sultanbekov, F.F. {\it Finite concrete
logics: their structure and measures on them}, Int. J.
Theor.Phys., 37(1998), 147--153.

9. Ovchinnikov, P.,{\it Measures on finite concrete logics},
Proc. Amer.Math. Soc.,127 (1999), 1957--1966.

10. Ovchinnikov, P.G., {\it On homogeneous finite Greechie logics
admiting a two-valued state,} In: Teor. funktsii, prilozh. i sm.
vopr., Kazan State University, Kazan, (1999), 167--168
(Russian).

11. Ptak, P.,{\it Exotic logics}, Coll. Math. 56 (1987), 1--7.

12. Rogalewicz, V. {\it A remark on $\lambda$-regular orthomodular
lattices}, Aplikace Mat., 34 (1989), 449--452.

13. Sultanbekov, F.F. {\it Signed measures and automorphisms of a
class of finite concrete logics}, Konstr. Teoriy funktsii i Funk
Analiz (Kazan), 8(1992), 57--68 (in Russian).

\vskip 0.1 cm
Mail: Department of Mathematics and Mechanics, Kazan State University,
Kremlevskaja Street, 18. Kazan, Tatarstan, Russia, 420008.

E-mail: Foat.Sultanbekov@ksu.ru

\end{document}